\documentclass[11pt]{article}
\usepackage[colorlinks=true,
urlcolor=blue,linkcolor=blue,citecolor=blue]{hyperref}
\usepackage{amsmath,amsfonts,amssymb,amsthm}
\usepackage{mathrsfs}
\usepackage{color}
\usepackage{amscd}
\newtheorem{theorem}{Theorem}[section]
\newtheorem{lemma}[theorem]{Lemma}

\newtheorem{proposition}[theorem]{Proposition}

\theoremstyle{remark}

\numberwithin{equation}{section}

\def\bxi{ {\boldsymbol\xi} }

\title{Convergence to SPDE of the Schr{\"o}dinger equation with large, random potential}
\author{Ningyao Zhang, Guillaume Bal}
\begin{document}
\maketitle

\begin{abstract}
	We study the asymptotic behavior of solutions to the  Schr{\"o}dinger equation with large-amplitude, highly oscillatory, random potential. In dimension $d<\mathfrak{m}$, where $\mathfrak{m}$ is the order of the leading operator in the Schr\"odinger equation, we construct the heterogeneous solution by using a Duhamel expansion and prove that it converges in distribution, as the correlation length $\varepsilon$ goes to $0$, to the solution of a stochastic differential equation, whose solution is represented as a sum of iterated Stratonovich integral, over the space $C([0,+\infty),\mathcal{S}')$. The uniqueness of the limiting solution in a dense space of $L^2(\Omega\times\mathbb{R}^d)$ is shown by verifying the property of conservation of mass for the Schr\"odinger equation. In dimension $d>\mathfrak{m}$, the solution to the Schr{\"o}dinger equation is shown to converge in $L^2(\Omega\times\mathbb{R}^d)$ to a deterministic Schr{\"o}dinger solution in \cite{ZB-12}.   
\end{abstract}

\section{Introduction}
We consider the following Schr{\"o}dinger equation in dimension $d<\mathfrak{m}$:
\begin{equation}
\label{eq1.1}
\left\{
\begin{aligned}
\Big(i\frac{\partial}{\partial t}+\big(P(D)-\frac{1}{\varepsilon^{d/2}}q(\frac{x}{\varepsilon})\big)\Big)u_\varepsilon(t,x)&=0, \qquad t>0,\, x\in \mathbb{R}^d\\
u_\varepsilon(0,x)&=u_0(x), \qquad x\in \mathbb{R}^d,
\end{aligned}
\right.
\end{equation}
where $P(D)$ is the pseudo-differential operator with symbol $\hat{p}(\xi)=|\xi|^{\mathfrak{m}}$. Taking Fourier transform of both sides of \eqref{eq1.1}, we obtain
\begin{equation}
\label{eq1.2}
\left\{
\begin{aligned}
(i\frac{\partial}{\partial t}+\xi^{\mathfrak{m}})\hat{u}_{\varepsilon}&=\varepsilon^{-\frac{d}{2}}
\int \hat{q}(\zeta) \hat{u}_{\varepsilon}(t,\xi-\varepsilon^{-1}\zeta)d\zeta,\\
\hat{u}_{\varepsilon}(0,\xi)&=\hat{u}_0(\xi).
\end{aligned}
\right.
\end{equation}
We assume the Fourier transform of the covariance of the potential $\hat{R}(\xi)$ is bounded and continuous at $0$, and the initial condition satisfies $(1+|\xi|^{2\mathfrak{m}})|\hat{u}_0(\xi)|\leq C$ uniformly in $\xi\in\mathbb{R}^d$. \\

The main objective of this paper is to construct a solution to the above equation in $L^2(\Omega\times\mathbb{R}^d)$ uniformly in time on bounded intervals and to show that the solution converges in distribution as $\epsilon\rightarrow 0$ to the unique solution of the following stochastic partial differential equation (SPDE)
\begin{equation}
	\label{SPDE}
\left\{
\begin{aligned}
i\frac{\partial u}{\partial t}+P(D)u-\sigma u\circ \dot{W}&=0, \qquad t>0,\, x\in \mathbb{R}^d\\
u(0,x)&=u_0(x), \qquad x\in \mathbb{R}^d,
\end{aligned}
\right.
\end{equation}
where $\dot{W}$ denotes spatial white noise, $\circ$ denotes the Stratonovich product and $\sigma$ is defined as
\begin{equation}
\sigma^2:=(2\pi)^{d}\hat{R}(0)=(2\pi)^d\int_{\mathbb{R}^d}R(x)dx.
\end{equation}
To study this solution, we may take the Fourier transform of the equation:
\begin{equation}
\label{eq2.2}
\begin{aligned}
(i\frac{\partial}{\partial t}+\xi^{\mathfrak{m}})\hat{u}(t,\xi)&=\sigma(2\pi)^{-d}\int e^{-i\xi x}
u(s,x)\circ dW(x)\\
&=\sigma(2\pi)^{-d}\int e^{-i\xi x}\int e^{i\xi_1x}\hat{u}(s,\xi_1)d\xi_1\circ dW(x),
\end{aligned}
\end{equation}
with initial condition $\hat{u}(0,\xi)=\hat{u}_0(\xi)$. To look for a mild solution, we recast \eqref{eq2.2} as 
\begin{equation}
\label{eq2.3}
\begin{aligned}
	\hat{u}=i\sigma(2\pi)^{-d}\int \int_0^t e^{i\xi^{\mathfrak{m}}(t-s)} e^{-i\xi x} \int e^{i\xi_1 x} \hat{u}(s,\xi_1)d\xi_1 ds \circ dW(x)+e^{i t \xi^{\mathfrak{m}}}\hat{u}_0(\xi).
\end{aligned}
\end{equation}
Define formally the stochastic integral 
\begin{equation}
	\mathcal{H}\hat{u}(t,\xi)=(-i\sigma)(2\pi)^{-d}\int_0^t e^{i(t-s)\xi^{\mathfrak{m}}}
\int e^{-i\xi x}\int e^{i\xi_1 x} \hat{u}(s,\xi_1) d\xi_1ds\circ dW(x).
\end{equation}
We may rewrite \eqref{eq2.3} as
\begin{equation}
	\hat{u}(t,\xi)=e^{it\xi^{\mathfrak{m}}}\hat{u}_0(\xi)+\mathcal{H}\hat{u}(t,\xi).
\end{equation}
The mild solution to \eqref{SPDE} is thus defined as $u(t,x)=\mathcal{F}^{-1}\{\hat{u}(t,\xi)\}$, where $\mathcal{F}^{-1}$ denotes inverse Fourier transform. Suppose that $d<\mathfrak{m}$. The following result holds for the solution of \eqref{eq2.2}.
\begin{theorem}
	Suppose $d<\mathfrak{m}$. The series
\begin{equation}
\label{limiting}
\hat{u}(t,\xi):=\sum_{n\geq 0}\hat{u}^{(n)}(t,\xi)
\end{equation}
converges in the $L^2(\Omega\times\mathbb{R}^d)$ sense for each $t\geq 0$ and $\xi\in \mathbb{R}^d$ and is the unique solution to \eqref{eq2.2} in the space $M$ dense in $L^2(\Omega\times\mathbb{R}^d)$, which is defined in Section \ref{Uniqueness}, where
\begin{equation}
\label{thm1}
\begin{aligned}
\hat{u}^{(n)}=&(-i\sigma)^{n}(2\pi)^{-nd} e^t\int_0^{\infty}\cdots\int_0^{\infty}\int e^{i\beta t}d\beta \int\prod_{k=1}^n d\xi_k\\
&\left\{\prod_{k=0}^{n}\left[1-i\left(\left|\xi-\sum_{j=1}^{k}\xi_j\right|^{\mathfrak{m}}-\beta\right)\right]\right\}^{-1}
\prod_{j=1}^{n} e^{-i\xi_j x_j} \hat{u}_0(\xi-\sum_{j=1}^n \xi_j)\circ\prod_{j=1}^{n} dW(\xi_j).
\end{aligned}
\end{equation}
\end{theorem}
The following theorem shows the weak convergence of $u_{\varepsilon}(t,x)$ to $u(t,x)$ for any $t>0$ and $x\in \mathbb{R}^d$.
\begin{theorem}
	Suppose that $d<\mathfrak{m}$. For any integers $r\geq 1$, $m_1,\cdots,m_r\geq 0$ and $\xi^{(1)},\cdots,\xi^{(r)}\in \mathbb{R}^d$, $t_1, \cdots, t_r \geq 0$ we have the convergence of moments
\begin{equation}
\label{thm2}
\lim_{\varepsilon\rightarrow 0_+}\mathbb{E}\{[\hat{u}_{\varepsilon}(t_1,\xi^{(1)})]^{m_1}\cdots [\hat{u}_{\varepsilon}(t_r,\xi^{(r)})]^{m_r}\}=\mathbb{E}\{[\hat{u}(t_1,\xi^{(1)})]^{m_1}\cdots[\hat{u}(t_r,\xi^{(r)})]^{m_r}\}.
\end{equation}
The finite dimensional distribution of $\hat{u}(t,\xi)$ is uniquely determined by its moments of all orders. Moreover, the family of processes $\{\hat{u}_{\varepsilon}(t,\cdot),t\geq 0\}$ is tight, as $\varepsilon\rightarrow 0_+$, over $C([0,+\infty);\mathcal{S}'(\mathbb{R}^d))$. The process $\{\hat{u}_{\varepsilon}(t,\cdot),t\geq 0\}$ converges in law over $C([0,+\infty);\mathcal{S}'(\mathbb{R}^d))$, as $\varepsilon\rightarrow 0_+$, to $\{\hat{u}(t,\cdot),t\geq 0\}$. Also, we have that in the spatial domain, the process $\{u_{\varepsilon}(t,\cdot),t\geq 0\}$ converges in law to $\{u(t,\cdot),t\geq 0\}$. 
\end{theorem}

The rest of the paper is structured as follows. Section 2 gives the formal Duhamel solutions for both the multi-scale Schr\"odinger equation and the limiting SPDE in the Fourier domain. Section 3 demonstrates the first order moment convergence of Duhamel solutions $\hat{u}_{\varepsilon}$. Section 4 proves that $\hat{u}$ defined as the Duhamel expansion of the limiting equation is the space of in $L^2(\Omega\times\mathbb{R}^d)$. Section 5 generalizes the first order moment convergence proved in Section 3 to arbitrary orders. Section 6 follows the approach as in \cite{KN-PA-10} to show the weak convergence of $\{\hat{u}_{\varepsilon}(t,\xi)\}$ in $C([0,+\infty),\mathcal{S}')$ to $\hat{u}(t,\xi)$ follows from tightness and convergence in finite dimensional distribution.  

The treatment of Schr{\"o}dinger equation \eqref{eq1.1} in dimension $d>\mathfrak{m}$ with potential with corresponding amplitude is presented in \cite{ZB-12}. It is then shown that $u_{\varepsilon}$ converges in $L^2(\Omega\times\mathbb{R}^d)$ to the solution of a homogenized equation. The asymptotic theory of solution to parabolic equation with large potential is presented in \cite{B-CMP-09, B-MMS-10}.     

\section{Duhamel Expansion}
Iteratively using Duhamel's formula we obtain
\begin{equation}
\label{ue1}
\begin{aligned}
\hat{u}_{\varepsilon}(t,\xi)&=e^{i\xi^{\mathfrak{m}}t}\hat{u}_0(\xi)-i\varepsilon^{-\frac{d}{2}}\int_0^t
\int e^{i\xi^{\mathfrak{m}}(t-s)}\hat{q}(\zeta)\hat{u}_{\varepsilon}(s,\xi-\varepsilon^{-1}\zeta)ds d\zeta\\
&=\sum_{n=0}^{+\infty} \hat{u}_{\varepsilon}^{(n)}(t,\xi),
\end{aligned}
\end{equation}
where
\begin{equation}
\begin{aligned}
\label{uen1}
\hat{u}_{\varepsilon}^{(n)}(t,\xi)=&(-i)^{n}\varepsilon^{-\frac{nd}{2}}\int\cdots\int_{\Delta_n(t)} ds_1\cdots ds_n\int\cdots\int \prod_{k=1}^{n+1} e^{i(s_{k-1}-s_{k})\left|\xi-\varepsilon^{-1}\sum_{j=1}^{k-1}\xi_j\right|^{\mathfrak{m}}}\\
&\prod_{k=1}^{n} \hat{q}(\xi_k)d\xi_k \hat{u}_0(\xi-\varepsilon^{-1}\sum_{j=1}^{n}\xi_j).
\end{aligned}
\end{equation}
Here, we introduce the notation $\sum_{j=1}^{0}\xi_j:=0$, and $\Delta_n(t):=[t\geq s_1\geq \cdots \geq s_n\geq 0]$. Let $\tilde{\Delta}_n(t):=[\sum_{j=1}^n \tau_j\leq t,\tau_j\geq 0]$. Changing variables $s_j:=\sum_{i=j}^n \tau_i$ and denoting $\tau_0:=t-\sum_{i=1}^n \tau_i$ we can rewrite \eqref{uen1} in the form
\begin{equation}
\begin{aligned}
\hat{u}_{\varepsilon}^{(n)}(t,\xi)=&(-i)^n \varepsilon^{-\frac{nd}{2}} \int\cdots\int_{\tilde{\Delta}_{n(t)}} d\tau_1\cdots d\tau_n \int\cdots\int \prod_{k=1}^{n} \hat{q}(\xi_k)d\xi_k \\
&\times\prod_{k=1}^{n+1} e^{i\tau_{k-1}\left|\xi-\varepsilon^{-1}\sum_{j=1}^{k-1} \xi_j\right|^{\mathfrak{m}}} \hat{u}_0(\xi-\varepsilon^{-1}\sum_{j=1}^{n}\xi_j)\\
=&(-i)^{n} \varepsilon^{-\frac{nd}{2}}\int_0^{+\infty}\cdots\int_{0}^{+\infty}\int\cdots\int d\tau_0\cdots d\tau_n \delta(t-\sum_{j=0}^n\tau_j)\prod_{k=1}^n \hat{q}(\xi_k)d\xi_k\\
&\times\prod_{k=1}^{n+1} e^{i\tau_{k-1}\left|\xi-\varepsilon^{-1}\sum_{j=1}^{k-1} \xi_j\right|^{\mathfrak{m}}}.
\end{aligned}
\end{equation}
Using $\delta(t)=\int e^{i\beta t} d\beta$, we obtain for any $\eta>0$
\begin{equation}
\begin{aligned}
\hat{u}_{\varepsilon}^{(n)}(t,\xi)=&(-i)^n \varepsilon^{-\frac{nd}{2}} e^{\eta t} \int_0^{+\infty}\cdots\int_{0}^{+\infty}\int\cdots\int d\tau_0\cdots d\tau_n \prod_{k=1}^n \hat{q}(\xi_k) d\xi_k\\
&\times e^{i\beta(t-\sum_{j=0}^n \tau_j)}\prod_{k=1}^{n+1} e^{i\tau_{k-1}
\left|\xi-\varepsilon^{-1}\sum_{j=1}^{k-1}\xi_j\right|^{\mathfrak{m}}} \hat{u}_0(\xi-\varepsilon^{-1}\sum_{j=1}^n \xi_j).
\end{aligned}
\end{equation}
Integrating out all $\tau_j$ and choosing $\eta=1$ we get
\begin{equation}
\label{uen2}
\begin{aligned}
\hat{u}_{\varepsilon}^{(n)}(t,\xi)=&(-i)^n \varepsilon^{-\frac{nd}{2}} e^{\eta t}\int_{\mathbb{R}}e^{i\beta t} d\beta \int\cdots\int \prod_{k=1}^n \hat{q}(\xi_k)d\xi_k\\
&\times \left\{\prod_{k=0}^n\left[1-i\left(\left|\xi-\varepsilon^{-1}\sum_{j=1}^k \xi_j\right|^{\mathfrak{m}}-\beta\right)\right]\right\}^{-1} \hat{u}_0(\xi-\varepsilon^{-1}\sum_{j=1}^n \xi_j).
\end{aligned}
\end{equation}

We now come to the analysis of the limiting equation. By Duhamel's principle, the solution to \eqref{eq2.2} formally satisfies the equation
\begin{equation}
\label{sol2.2}
\begin{aligned}
\hat{u}(t,\xi)=&e^{i\xi^{\mathfrak{m}}t}\hat{u}_0(\xi)+(-i\sigma)(2\pi)^{-d}\int_0^t e^{i(t-s)\xi^{\mathfrak{m}}}
\int e^{-i\xi x} u(s,x) ds\circ dW(x)\\
=&e^{i\xi^{\mathfrak{m}}t}\hat{u}_0(\xi)+(-i\sigma)(2\pi)^{-d}\int_0^t e^{i(t-s)\xi^{\mathfrak{m}}}
\int e^{-i\xi x}\int e^{i\xi_1 x} \hat{u}(s,\xi_1) d\xi_1ds\circ dW(x).
\end{aligned}
\end{equation}
Integrating \eqref{sol2.2} iteratively, we obtain formally the Duhamel expansion for \eqref{eq2.2}
\begin{equation}
	\label{Duhamel}
\hat{u}(t,\xi)=\sum_{n=0}^{\infty} \hat{u}^{(n)},
\end{equation}
where
\begin{equation}
	\hat{u}^{(0)}=e^{it\xi^{\mathfrak{m}}}\hat{u}_0(\xi), \qquad \text{and} \qquad \hat{u}^{(n)}=\mathcal{H}^n\hat{u}^{(0)},
\end{equation}
for $n=0,1,\cdots$, or more explicitly,
\begin{equation}
\begin{aligned}
	\label{stratonovich}
\hat{u}^{(n)}=&(-i\sigma)^{n}(2\pi)^{-nd} \int\cdots\int_{\Delta_n(t)}\int\int\prod_{k=1}^{n+1}e^{i(s_{k-1}-s_k)\xi_{k-1}^{\mathfrak{m}}}\\
&\prod_{j=1}^n e^{-i(\xi_j-\xi_{j-1})x_j} \hat{u}_0(\xi_n) \prod_{k=1}^{n}d\xi_k \prod_{k=1}^{n}ds_k \circ\prod_{j=1}^{n} dW(x_j).
\end{aligned}
\end{equation}
Here, we define $s_0=t$, $s_{n+1}=0$, and $\xi_0:=\xi$. 
At this point neither the iterative Stratonovich integral of $\hat{u}^{(n)}$ for $n=1,2,\cdots$ in \eqref{stratonovich} nor the sum of \eqref{Duhamel} are well defined. We give a justification for these expressions in section 4. \\
Using the change of variables $\xi_k\rightarrow \xi_{k-1}-\xi_{k}$, and applying the same type of transform as in \eqref{uen2}, we obtain
\begin{equation}
\label{un2}
\begin{aligned}
	\hat{u}^{(n)}=&(-i\sigma)^{n}(2\pi)^{-nd} e^t\int_{\mathbb{R}^d}\int\cdots\int e^{i\beta t}d\beta \int\prod_{k=1}^n d\xi_k\\
&\left\{\prod_{k=0}^{n}\left[1-i\left(\left|\xi-\sum_{j=1}^{k}\xi_j\right|^{\mathfrak{m}}-\beta\right)\right]\right\}^{-1}
\prod_{j=1}^{n} e^{-i\xi_j x_j} \hat{u}_0(\xi-\sum_{j=1}^n \xi_j)\circ\prod_{j=1}^{n} dW(x_j),
\end{aligned}
\end{equation}
which is the same as \eqref{thm1}. 

\section{Convergence of the first moments}
We shall prove the convergence of moments as stated in \eqref{thm2}. For simplicity, we shall first consider the case when $n=m=1$ and show that the limit
\begin{equation}
\lim_{\varepsilon\rightarrow 0_+}\mathbb{E}\hat{u}_{\varepsilon}(t,\xi)
\end{equation}
exists. The proof of convergence for general moments is given in Section \ref{section:general}.  

Taking expectation of both sides of \eqref{ue1}, we obtain
\begin{equation}
\mathbb{E}\hat{u}_{\varepsilon}(t,\xi)=\sum_{n\geq 0} \mathbb{E}\hat{u}_{\varepsilon}^{(2n)}(t,\xi).
\end{equation}
This is because expectation of product of odd number of Gaussian random variables is $0$. The expectation of product of even number of Gaussian random variables is given as a sum of products of the expectation of pairs of variables, where the summation runs over all possible pairs.The contribution of products of potentials can thus be represented by
\begin{equation}
\mathbb{E}\{\prod_{k=1}^{2n}\hat{q}(\xi_k)\}=\sum_{\pi}\prod_{(ef)\in \pi} \hat{R}(\xi_e)\delta(\xi_e+\xi_f),
\end{equation}
where $(ef)$ denotes pair of indices, $\pi$ denotes a pairing of the $2n$ indices, and the summation is over all possible pairings.\\
Adding up all the delta functions gives $\sum_{k=1}^{2n}\xi_k=0$. We can therefore write
\begin{equation}
\begin{aligned}
\mathbb{E}\hat{u}_{\varepsilon}^{(2n)}(t,\xi)=&(-1)^n \varepsilon^{{nd}} e^{t} \hat{u}_0(\xi) \int_{\mathbb{R}} e^{i\beta t}d\beta \int\cdots\int\\
&\times \left\{\prod_{k=0}^{2n} \left[1-i\left(\left|\xi-\varepsilon^{-1}\sum_{j=1}^k \xi_j\right|^{\mathfrak{m}}-\beta\right)\right]\right\}^{-1}\mathbb{E}\{\prod_{k=1}^{2n}\hat{q}(\xi_k)d\xi_k\}\\
=&(-1)^n \varepsilon^{{nd}} e^{t} \hat{u}_0(\xi) \sum_{\pi} \int_{\mathbb{R}} e^{i\beta t}d\beta \int\cdots\int\\
&\times \left\{\prod_{k=0}^{2n} \left[1-i\left(\left|\xi-\varepsilon^{-1}\sum_{j=1}^k \xi_j\right|^{\mathfrak{m}}-\beta\right)\right]\right\}^{-1}\prod_{(ef)\in \pi} \hat{R}(\xi_e)\delta(\xi_e+\xi_f) d\xi_e d\xi_f.\\
\end{aligned}
\end{equation}
The above summation extends over all possible pairings $\pi$ made over vertices $\{1,\cdots,2n\}$. By changing variables $\xi_k:=\varepsilon^{-1}\xi_k$, we obtain
\begin{equation}
\label{ue2n}
\begin{aligned}
\mathbb{E}\hat{u}_{\varepsilon}^{(2n)}(t,\xi)=&(-1)^n e^{t} \hat{u}_0(\xi) \sum_{\pi} \int_{\mathbb{R}} e^{i\beta t}d\beta \int\cdots\int\\
&\times \left\{\prod_{k=0}^{2n} \left[1-i\left(\left|\xi-\sum_{j=1}^k \xi_j\right|^{\mathfrak{m}}-\beta\right)\right]\right\}^{-1}\prod_{(ef)\in \pi} \hat{R}(\xi_e)\delta(\xi_e+\xi_f) d\xi_e d\xi_f.
\end{aligned}
\end{equation}
In what follows, we show that
\begin{equation}
	\label{2nbound}
	\sup|\mathbb{E}\hat{u}_{\varepsilon}^{(2n)}(t,\xi)|\leq a_n,
\end{equation}
where $\sum_n a_n<+\infty$. Hence,
\begin{equation}
\label{1stlimit}
\lim_{\varepsilon\rightarrow 0_+} \mathbb{E} \hat{u}_{\varepsilon}(t,\xi)=\sum_{n=0}^{+\infty}
\lim_{\varepsilon\rightarrow 0_+} \mathbb{E} \hat{u}_{\varepsilon}^{(2n)}(t,\xi).
\end{equation}
Denote by $\mathcal{L}(\pi)$ the set of all left vertices of a given pairing $\pi$. Let
\begin{equation}
A_k:=\left|\xi-\sum_{j=1}^k\xi_j\right|^{\mathfrak{m}}.
\end{equation}
Define
\begin{equation}
\hat{e}_A^{(\rho)}=[1-i(A-\beta)]^{-1}.
\end{equation}
Using contour integration method, we are able to show that
\begin{equation}
{e}_A^{(\rho)}=\left\{
\begin{array}{ll}
e^{iAt}t^{\rho-1}e^{-t} & t>0.\\
0 & t<0
\end{array}
\right.
\end{equation}
The following estimate is then derived:
\begin{lemma}
Suppose that $\rho>0$. There exists a constant $C>0$ such that for an arbitrary $n\geq 1$ we have
\begin{equation}
\left|\int e^{i\beta t}\left\{\prod_{k=1}^{n}[1-i(A_k-\beta)]\right\}^{-\rho}d\beta\right|\leq \frac{C^n t^{n\rho-1} e^{-t}}
{[(n-1)!]^{\rho}}.
\end{equation}
\end{lemma}
\noindent Readers are referred to \cite{KN-PA-10} for the proof of this lemma.\\
With such notation, we may rewrite \eqref{ue2n} in the form
\begin{equation}
\begin{aligned}
\mathbb{E} \hat{u}_{\varepsilon}^{(2n)} (t,\xi) = &(-1)^n e^t \hat{u}_0(\xi) \sum_{\pi} \int\cdots\int
\prod_{k\notin \mathcal{L}(\pi), k\neq 2n}*{e}_{A_k}*F(t,\bxi;\pi)\\
&\times\prod_{(ef)\in \pi} \hat{R}(\varepsilon\xi_e)\delta(\xi_e+\xi_f)d\xi_e d\xi_f,
\end{aligned}
\end{equation}
where
\begin{equation}
F(t,\bxi;\pi):=\int e^{i\beta t} (1-i(\xi^{\mathfrak{m}}-\beta))^{-2}\prod_{k\in \mathcal{L}(\pi)}
\left[1-i\left(\left|\xi-\sum_{j=1}^k \xi_j\right|^{\mathfrak{m}}-\beta\right)\right]^{-1}d \beta
\end{equation}
and $\prod_{k\notin \mathcal{L}(\pi), k\neq 2n}*\hat{e}_{A_k}$ denotes the convolution of all ${e}_{A_k}$'s, where $k$ is the right endpoint in the giving pairing $\pi$, except $k=2n$.
\begin{lemma}
Suppose $\rho\in(\frac{d}{\mathfrak{m}},1)$. Then
\begin{equation}
\sup_{\beta\in\mathbb{R},\omega\in\mathbb{R}^d} \int_{\mathbb{R}^d}
|-\beta+|\xi-\omega|^{\mathfrak{m}}+i|^{-\rho}d\xi<+\infty.
\end{equation}
\end{lemma}
\begin{proof}
We may first shift $\xi$ to get rid of $\omega$ and then change
variables to polar coordinates
\begin{equation}
\label{lemma3.2.1}
\int_{\mathbb{R}^d}
|-\beta+|\xi-\omega|^{\mathfrak{m}}+i|^{-\rho}\leq
\Omega_d\int_0^{+\infty}
|-\beta+\xi^{\mathfrak{m}}+i|^{-\rho}\xi^{d-1} d|\xi|,
\end{equation}
where $\Omega_d$ is the area of the unit sphere in $\mathbb{R}^d$.
Let $Q=\xi^{\mathfrak{m}}$. The integral on the right hand side of
\eqref{lemma3.2.1} may be rewritten as
\begin{equation}
\int_0^{+\infty} |-\beta+Q+i|^{-\rho}Q^{d/{\mathfrak{m}}-1} dQ.
\end{equation}
Without loss of generality, we assume $\beta\geq0$ and the integral
above can then be written as the sum of three integrals $I$, $II$,
$III$ according to whether $\xi$ belongs to $(0,\beta/2)$,
$(\beta/2,2\beta)$, or $(2\beta,+\infty)$. We have
\begin{equation}
I\leq
[1+(\frac{\beta}{2})^2]^{-\frac{\rho}{2}}\int_0^{\beta/2}Q^{\frac{d}{\mathfrak{m}}-1}dQ
\leq
\frac{\mathfrak{m}}{d}[1+(\frac{\beta}{2})^2]^{-\frac{\rho}{2}}(\frac{\beta}{2})^{\frac{d}{\mathfrak{m}}}\leq
\frac{\mathfrak{m}}{d}.
\end{equation}
The second integral can also be estimated as follows
\begin{equation}
II\leq
(2\beta)^{\frac{d}{\mathfrak{m}}-1}\int_{\frac{\beta}{2}}^{2\beta}|Q-\beta+i|^{-\rho}dQ.
\end{equation}
If $\beta\leq 1$ we have
\begin{equation}
II\leq 3\times 2^{\frac{d}{\mathfrak{m}}-2}.
\end{equation}
If $\beta>1$, we estimate
\begin{equation}
II\leq
(2\beta)^{\frac{d}{\mathfrak{m}}-1}\int_{\frac{\beta}{2}}^{2\beta}|Q-\beta|^{-\rho}dQ\leq
[1+(\frac{1}{2})^{1-\rho}]2^{\frac{d}{\mathfrak{m}}-1}.
\end{equation}
The third integral is estimated
\begin{equation}
III\leq
\frac{\mathfrak{m}}{d}(\beta^2+1)^{-\frac{\rho}{2}}\beta^{\frac{d}{\mathfrak{m}}}\leq\frac{\mathfrak{m}}{d}.
\end{equation}This concludes the proof of the lemma.
\end{proof}
Let
\begin{equation}
	\label{F1}
F_1(t,\bxi;\pi):=\int_{\mathbb{R}}\frac{e^{i\beta t} d\beta}{[1-i(\xi^{\mathfrak{m}}-\beta)]^2}
\left\{\prod_{k\in\mathcal{L}(\pi)}\left[1-i\left(\left|\xi-\sum_{j=1}^{k}\xi_j\right|^{\mathfrak{m}}-\beta\right)\right]^{-\rho}\right\}
\end{equation}
and
\begin{equation}
F_2(t,\bxi;\pi):=\int_{\mathbb{R}}e^{i\beta t}d\beta
\left\{\prod_{k\in\mathcal{L}(\pi)}\left[1-i\left(\left|\xi-\sum_{j=1}^{k}\xi_j\right|^{\mathfrak{m}}-\beta\right)\right]^{-(1-\rho)}\right\}.
\end{equation}
Since
\begin{equation}
F_1(t,\bxi;\pi)=e^{(\rho)}_{B_1(\bxi)}*\cdots *e^{(\rho)}_{B_{n+2}(\bxi)} (t)
\end{equation}
for some $B_j(\bxi)$, we get $F_1(t,\bxi;\pi)=0$ for $t<0$. Likewise, $F_2(t,\bxi;\pi)=0$ for $t<0$. We can therefore write
\begin{equation}
F(t,\bxi;\pi)=\int_0^t F_1(t-s,\bxi;\pi) F_2(s,\bxi;\pi) ds\geq 0.
\end{equation}
Observe that by Lemma 3.1
\begin{equation}
F_2(t,\bxi;\pi)\leq \frac{C^n t^{n(1-\rho)-1}}{[(n-1)!]^{1-\rho}}
\end{equation}
and
\begin{equation}
\prod_{k\notin \mathcal{L}(\pi),k\neq 2n} * e_{A_k}(t)\leq \frac{C^n t^{n-1}}{(n-1)!}
\end{equation}
for $t>0$.
On the other hand, we have obviously
\begin{equation}
	|F_1(t,\bxi;\pi)|\leq C\int_{\mathbb{R}} \frac{d\beta}{1+|\xi^{\mathfrak{m}}-\beta|^2}
\left\{\prod_{k\in\mathcal{L}(\pi)}\left|1-i\left(\left|\xi-\sum_{j=1}^{k}\xi_j\right|^{\mathfrak{m}}-\beta\right)\right|^{-\rho}\right\}
\end{equation}
for some constant $C>0$. As a result, we obtain
\begin{equation}
\begin{aligned}
	|\mathbb{E}\hat{u}_{\varepsilon}^{(2n)}(t,\xi)|\leq &C e^t |\hat{u}_0(\xi)| \sum_{\pi} \int_0^t \frac{(t-s)^{n-1}ds}{(n-1)!}\\
&\times \int\cdots\int F(s,\bxi;\pi) \prod_{(kl)\in\pi} \hat{R}(\varepsilon\xi_e) \delta(\xi_e+\xi_f) d\xi_e d\xi_f\\
&\leq \frac{C^n e^t |\hat{u}_0(\xi)|}{[(n-1)!]^{2-\rho}} \int_0^t (t-s)^n s^{n(1-\rho)-1} ds \sum_{\pi} \int_{\mathbb{R}} \frac{d\beta}{1+|\xi^{\mathfrak{m}}-\beta|^2} \int\cdots\int\\
&\times \prod_{(ef)\in\pi}\hat{R}(\varepsilon\xi_e)\delta(\xi_e+\xi_f)d\xi_e d\xi_f
\left\{\prod_{k\in\mathcal{L}(\pi)}\left|1-i\left(\left|\xi-\sum_{j=1}^{k}\xi_j\right|^{\mathfrak{m}}-\beta\right)\right|^{-\rho}\right\}.
\end{aligned}
\end{equation}
Here, we just need half of the order of decay of the initial condition $\hat{u}_0(\xi)$ as assumed to have
\begin{equation}
	\frac{|\hat{u}_0(\xi)|}{(1+|\xi^{\mathfrak{m}}-\beta|^2)^{1/2}}\leq \frac{C}{(1+\xi^{2\mathfrak{m}})^{1/2}(1+|\xi^{\mathfrak{m}}-\beta|^2)^{1/2}}\leq \frac{1}{(1+\beta^2)^{1/2}}.
\end{equation}
The right hand side of the above inequality is then used for estimating the integration in $\beta$ by Lemma 5.2 in \cite{ZB-12}:
\begin{equation}
	\int_{-\infty}^{+\infty} \frac{d\beta}{(1+\beta^2)^{1/2}(1+|\xi^{\mathfrak{m}}-\beta|^2)^{1/2}}\leq C\frac{1+\log_{+}|\xi|}{(1+\xi^{2\mathfrak{m}})^{1/2}},
\end{equation}
where $\log_{+} |\xi|:=max(0,\log |\xi|)$.\\
Using Lemma 3.2 we conclude therefore that
\begin{equation}
\label{1stue}
|\mathbb{E} \hat{u}_{\varepsilon}^{(2n)} (t,\xi)|\leq (2n-1)!!\frac{C^n t^{n(2-\rho)}e^t}{(n!)^{2-\rho}}\frac{1+\log_{+}|\xi|}{(1+\xi^{2\mathfrak{m}})^{1/2}}
\end{equation}
and \eqref{2nbound} follows by letting $a_{2n}=(2n-1)!!\frac{C^n t^{n(2-\rho)}e^t}{(n!)^{2-\rho}}$. We obtain
\begin{equation}
\lim_{\varepsilon\rightarrow 0_{+}} \mathbb{E}\hat{u}_{\varepsilon}(t,\xi)=\sum_{n=0}^{+\infty} \bar{u}^{(2n)} (t,\xi),
\end{equation}
where
\begin{equation}
\label{1stu}
\begin{aligned}
\bar{u}^{(2n)}(t,\xi):=&(-\hat{R}(0))^n e^t \hat{u}_0(\xi) \sum_{\pi} \int_{\mathbb{R}} e^{i\beta t}d\beta \int\cdots\int \\
&\times \left\{\prod_{k=0}^{2n}\left[1-i\left(\left|\xi-\sum_{j=1}^{k}\xi_j\right|^{\mathfrak{m}}-\beta\right)\right]\right\}^{-1}
\prod_{(ef)\in\pi}\delta(\xi_e+\xi_f)d\xi_e d\xi_f.
\end{aligned}
\end{equation}
In what follows, we show that
\begin{equation}
\label{equality}
\bar{u}^{(2n)}(t,\xi)=\mathbb{E}\hat{u}^{(2n)}(t,\xi).
\end{equation}
Expectation of multiple Stratonovich integrals is defined:
\begin{equation}
\begin{aligned}
\mathbb{E}\{\prod_{j=1}^{2n}dW(x_j)\}=\sum_{\pi}\prod_{(ef)\in\pi}\delta(x_e-x_f)dx_e dx_f.
\end{aligned}
\end{equation}
Upon integrating in all variables $x$, the first moment of $\hat{u}^{(2n)}$ can therefore be written as
\begin{equation}
\label{1stE}
\begin{aligned}
	\mathbb{E}\hat{u}^{(2n)}=&(-i\sigma)^{2n}(2\pi)^{-2nd} e^t\sum_{\pi}\int_0^{\infty}\cdots\int_0^{\infty}\int e^{i\beta t}d\beta \int\prod_{k=1}^{2n} d\xi_k\\
&\left\{\prod_{k=0}^{2n}\left[1-i\left(\left|\xi-\sum_{j=1}^{k}\xi_j\right|^{\mathfrak{m}}-\beta\right)\right]\right\}^{-1}
\hat{u}_0(\xi)\prod_{(ef)\in\pi}\delta(\xi_e+\xi_f),
\end{aligned}
\end{equation}
which gives exactly \eqref{equality}. 
\section{$L^2$ convergence of the SPDE solution}
We now come to the series \eqref{limiting}, which we claim is the solution to \eqref{eq2.3}. The readers are referred to \cite{B-CMP-09} for more complete description of the relationship between Stratonovich and It{\^o} integrals and for additional details on the theory presented in this section.

We first prove that $\hat{u}(t,\cdot)$ in \eqref{limiting} as a series is well-defined in the space of $L^2(\Omega\times\mathbb{R}^d)$. Denote $\hat{u}^{(n)}=\mathcal{I}_n(f_n)$, where $\mathcal{I}_n$ denotes the n-th order iterated Stratonovich integral, and $f_n$ is a n-parameter function, i.e.
\begin{equation}
	\mathcal{I}_n(f_n) = \int_{\mathbb{R}^{nd}} f_n(x_1,\dots,x_n) dW(x_1)\dots dW(x_n).
\end{equation}
By the definition of the $L^2$ norm of multiple Stratonovich integral, we have
\begin{equation}
	\begin{aligned}
	\mathbb{E}|\hat{u}^{(n)}|^2&=\mathbb{E}\{\mathcal{I}_{2n}(f_n\otimes \bar{f}_n)\}\\
	&=\mathbb{E}\int_{\mathbb{R}^{2nd}}f_n(x_1,\dots,x_n) \bar{f}_n(y_1,\dots,y_n)dW(x_1)\dots dW(x_n) dW(y_1)\dots dW(y_n).
\end{aligned}
\end{equation}
Taking the same approach as in the calculation of the first moment in \eqref{1stE}, we obtain
\begin{equation}
	\label{2ndmoment}
\begin{aligned}
	&\mathbb{E}|\hat{u}^{(n)}(t,\xi)|^2=(\sigma^{n}(2\pi)^{-d}e^t)^2\sum_{\pi_{2n}}\int\int \int_{\mathbb{R}^2} e^{i(\beta_1-\beta_2)t}
d\beta_1 d\beta_2 \hat{u}_0(\xi-\sum_{j=1}^{n} \xi^{(1)}_{j})\\
&\overline{\hat{u}}_0(\xi+\sum_{j=1}^{n}\xi_j^{(2)})\times\left\{\prod_{k=0}^{n}\left[1-i\left(\left|\xi-\sum_{j=1}^{k}\xi_j^{(1)}\right|^{\mathfrak{m}}+\beta_1\right)\right]
\prod_{k=0}^{n}\left[1+i\left(\left|\xi+\sum_{j=1}^{k}\xi_j^{(2)}\right|^{\mathfrak{m}}+\beta_2\right)\right]\right\}^{-1}\\
&\prod_{(ef)\in\pi}\delta(\xi_e+\xi_f)d\bxi^{(1)} d\bxi^{(2)} .
\end{aligned}
\end{equation}
We denote $d\bxi^{(l)}=d\xi_1^{(l)}\cdots d\xi_n^{(l)}, l=1,2$. $\xi_e$ and $\xi_f$ are paired arguments in the graph comprised of $2n$ arguments in total with their index $e$ and $f$.\\
The idea that we use to estimate the first moment applies here too. However, note that we have $\{k=n\}\in\mathcal{L}(\pi)$ for the crossing graphs defined in \cite{B-MMS-10}, and one of the terms $[1-i(\xi^{\mathfrak{m}}-\beta)]^{-1}$ in the function $F_1$ defined in \eqref{F1} now becomes $[1-i(|\xi-\sum_{j=1}^n\xi^{(1)}_j|^{\mathfrak{m}}-\beta)]^{-1}$. Therefore, we have to make the following adjustments in the proof. Using the smoothness condition for the initial condition, we have
\begin{equation}
	\left|\frac{\hat{u}_0(\xi-\sum_{j=1}^n\xi_j^{(1)})}{1-i(|\xi-\sum_{j=1}^{n}\xi_j^{(1)}|^\mathfrak{m}-\beta_1)}\right|\leq \frac{1}{(1+|\xi-\sum_{j=1}^n \xi_j^{(1)}|^{2\mathfrak{m}})^{1/2}}\frac{1}{(1+\beta_1^2)^{1/2}}.
\end{equation}
We use the term $[1+|\xi-\sum_{j=1}^n \xi_j^{(1)}|^{2\mathfrak{m}}]^{-1/2}$ for the integration in $\xi_n$, and the term $[1+\beta_1^2]^{-1/2}$ together with $[1-i(\xi^{\mathfrak{m}}+\beta_1)]^{-1}$ indexed by $\{k=0\}$ in the first product in \eqref{2ndmoment} 
\begin{equation}
	\int_{-\infty}^{+\infty} \frac{d\beta_1}{(1+\beta_1^2)^{1/2}|1-i(\xi^{\mathfrak{m}}+\beta_1)|}\leq C \frac{1+\log_{+}|\xi|}{(1+\xi^{2\mathfrak{m}})^{1/2}}.
\end{equation}
We also use the smoothness condition for the initial condition and obtain another $(1+\log_{+}|\xi|)/{(1+\xi^{2\mathfrak{m}})^{1/2}}$
from the integration in $\beta_2$.  
Finally, we obtain the estimate 
\begin{equation}
	\label{2ndbound}
	\mathbb{E}|\hat{u}^{(n)}(t,\xi)|^2\leq (2n-1)!!\frac{C^n t^{n(2-\rho)}e^t}{(n!)^{2-\rho}}\frac{1+\log_{+}|\xi|}{(1+\xi^{2\mathfrak{m}})^{1/2}}
\end{equation}
for any $\rho\in(\frac{d}{\mathfrak{m}},1)$. This implies that the iterated Stratonovich integral $\hat{u}^{(n)}$ is indeed well-defined. Integrating in $\xi$ and summing the above bound over $n$ gives
\begin{equation}
	\mathbb{E}\int|\hat{u}(t,\xi)|^2 d\xi = \sum_{n\geq 0} (\mathbb{E}\int|\hat{u}^{(n)}(t,\xi)|^2d\xi)^{1/2}<\infty
\end{equation}
and the $L^2(\Omega\times\mathbb{R}^d)$ convergence of \eqref{limiting} follows. In fact, by first multiplying $\xi^{\mathfrak{m}}$ with \eqref{2ndbound} and performing the integration and summation we can futher obtain
\begin{equation}
	\mathbb{E}\int\xi^{\mathfrak{m}}|\hat{u}(t,\xi)|^2 d\xi<\infty.
\end{equation}

\section{Uniqueness of the SPDE solution}\label{Uniqueness}
Let us now provide the rigorous definition for the operator $\mathcal{H}$. Suppose $f(t,\xi)$ is a sum of iterated Stratonovich integrals
\begin{equation}
	\label{expansion}
	f(t,\xi)=\sum_{n\geq 0} \mathcal{I}_n(f_n(t,\xi,\cdot)).
\end{equation}
We define
\begin{equation}
\mathcal{H}f(t,\xi)=\sum_{n \geq 1}\mathcal{I}_n( (\mathcal{H}f)_n(t,\xi,\cdot)),
\end{equation}
where
\begin{equation}
	(\mathcal{H}f)_{n+1}(t,\xi,x,x_1,\cdots,x_{n}):=(-i\sigma)(2\pi)^{-d}\int_0^t \int e^{i(t-s)\xi^{\mathfrak{m}}} e^{i(\xi_0-\xi)x} f_n(s,\xi_0,x_1,\cdots,x_n)d\xi_0 ds.
\end{equation}
It is checked that under this definition $\mathcal{H}\hat{u}=\sum_{n\geq 1}\hat{u}^{(n)}=\hat{u}-\hat{u}^{(0)}$. The Duhamel solution \eqref{limiting} is therefore a  solution to the equation \eqref{eq2.3}.

We can also define
\begin{equation}
	\mathcal{J}f(t,\xi)=\sum_{n\geq 0} \mathcal{I}_{n+1}( (\mathcal{J}f)_{n+1}(t,\xi,\cdot)),
\end{equation}
where
\begin{equation}
(\mathcal{J}f)_{n+1}(t,\xi,x,x_1,\cdots,x_{n}):=(-i\sigma)(2\pi)^{-d} e^{-i\xi x}\int  e^{i\xi_0 x} f_n(s,\xi_0,x_1,\cdots,x_n)d\xi_0.
\end{equation}
For the Duhamel solution \eqref{limiting}, we have 
\begin{equation}
	\begin{aligned}
	&\mathbb{E}|\mathcal{I}_{n+1}\left((\mathcal{J}\hat{u})_{n+1}(t,\xi)\right)|^2 = (\sigma^{n}(2\pi)^{-d}e^t)^2\sum_{\pi_{2n}}\int\int \int_{\mathbb{R}^2} e^{i(\beta_1-\beta_2)t}
d\beta_1 d\beta_2 \hat{u}_0(\xi-\sum_{j=1}^{n+1} \xi^{(1)}_{j})\\
&\times\overline{\hat{u}}_0(\xi+\sum_{j=1}^{n+1}\xi_j^{(2)})\left\{\prod_{k=1}^{n+1}\left[1-i\left(\left|\xi-\sum_{j=1}^{k}\xi_j^{(1)}\right|^{\mathfrak{m}}+\beta_1\right)\right]
\prod_{k=1}^{n+1}\left[1+i\left(\left|\xi+\sum_{j=1}^{k}\xi_j^{(2)}\right|^{\mathfrak{m}}+\beta_2\right)\right]\right\}^{-1}\\
&\prod_{(ef)\in\pi}\delta(\xi_e+\xi_f)d\bxi^{(1)} d\bxi^{(2)} .
        \end{aligned}
\end{equation}
Although it looks a little different from \eqref{2ndmoment}, it can be estimated in the same way as
\begin{equation}
(\mathbb{E}|\hat{u}^{(n)}|^2)^{1/2}<(2n-1)!!\frac{C^n t^{n(2-\rho)}e^t}{(n!)^{2-\rho}}\frac{1+\log_{+}|\xi|}{(1+\xi^{2\mathfrak{m}})^{1/2}},
\end{equation}
which uppon summation in $n$ implies that $\mathbb{E}|\hat{u}(t,\xi)|^2$ is uniformly bounded for all $\xi\in\mathbb{R}^d$. \\
Note that $\mathcal{H}\hat{u}(t,\xi) = (-i\sigma)(2\pi)^{-d}\int_0^t e^{i(t-s)\xi^{\mathfrak{m}}} (\mathcal{J}\hat{u})(s,\xi) ds$. $\hat{u}(t,\xi)$ is therefore a solution to the equation
\begin{equation}
	\label{eq2.2*}
	(i\frac{\partial}{\partial t}+\xi^{\mathfrak{m}})\hat{u}(t,\xi)=\sigma(2\pi)^{-d}\mathcal{J}\hat{u}(t,\xi).
\end{equation}
Now we prove that this equation preserves mass. Multiplying this equation by $\bar{\hat{u}}(t,\xi)$, and integrating in $\xi$ and over the probability space $\Omega$ gives
\begin{equation}
	\label{energy}
	\frac{i}{2}\frac{\partial }{\partial t}\mathbb{E}\int |\hat{u}|^2d\xi +\mathbb{E}\int\xi^{\mathfrak{m}}|\hat{u}|^2d\xi=\sigma(2\pi)^{-d} \mathbb{E}\int(\mathcal{J}\hat{u})\bar{\hat{u}}d\xi.
\end{equation}
The right hand side of this equation can be written out explicitly as
\begin{equation}
	\begin{aligned}
		\mathbb{E}\int (\mathcal{J}\hat{u}) \bar{\hat{u}} d\xi &= \int \sum_{n,m} \mathbb{E} \left\{\mathcal{I}_{n+1}\left(e^{-i\xi x}\int e^{i\xi_1 x}
		f_n(t,\xi_1)d\xi_1\right)\mathcal{I}_m \left(\bar{f}_m(t,\xi)\right)\right\}d\xi\\
		&=\int \sum_{n,m} \int \left( e^{-i\xi x} \int e^{i\xi_1 x} f_n(t,\xi_1) d\xi_1 \bar{f}_m(t,\xi) \sum_{\pi} \prod_{(ef)\in \pi} \delta(x_e-x_f) d{\bf x} \right) d\xi \\
		&=\sum_{n,m} \sum_{\pi} \int e^{i\xi_1 x} f_n(t,\xi_1) d\xi_1 \int e^{-i\xi x} \bar{f}_m(t,\xi) d\xi \prod_{(ef)\in \pi} \delta(x_e-x_f) d{\bf x} \\
		&=\sum_{n,m} \sum_{\pi} \int \check{f}_n(t,x) \bar{\check{f}}_m(t,x) \prod_{(ef)\in \pi} \delta(x_e-x_f) d{\bf x},
	\end{aligned}
\end{equation}
which is real-valued because of the symmetricity of this summation. Extracting the imaginary part from both sides of \eqref{energy} gives 
\begin{equation}
	\frac{\partial}{\partial t} \mathbb{E} \int |\hat{u}(t,\xi)|^2 d\xi = 0.
\end{equation}

Finally, we define the space $M$ in which the equation \eqref{eq2.2*} admits a unique solution. In light of the equation \eqref{energy}, $M$ consists of sum of iterated Stratonovich integrals $f(t,\xi) = \sum_{n\geq 0}\mathcal{I}_n(f_n(t,\xi,\cdot))$ such that

\begin{enumerate}
\item $f(t,\xi)\in L^2(\Omega\times\mathbb{R}^d)$,
\item $\mathcal{J}f \in L^2(\Omega)$,
\item $|\xi|^{\frac{\mathfrak{m}}{2}}f(t,\xi)\in L^2(\Omega\times\mathbb{R}^d)$.
\end{enumerate}

As a reminder, defining the sum of iterated Stratonovich integral $f$ as in \eqref{expansion}, we have that
\begin{equation}
	f = \sum_{n\geq 0}\mathcal{I}_n(f_n) = \sum_{m\geq 0}I_m(g_m),
\end{equation}
where
\begin{equation}\label{change}
	g_m(t,\xi,x) = \sum_{k\geq 0}\frac{(m+2k)!}{m!k!2^k}\int_{\mathbb{R}^{kd}}f_{m+2k}(t,\xi,{\bf x}_m,{\bf y}_k^{\otimes 2})d{\bf y}.
\end{equation}
Here, ${\bf y}\otimes {\bf y} \equiv (y,y)$, and 
\begin{equation}
I_m(g_m):=\int_{\mathbb{R}^{md}} g_m(t,\xi,x_1,\cdots,x_m)dW(x_1)\cdots dW(x_m)
\end{equation}
denotes the iterated It{\^o} integral. The $L^2$ norm of $f$ can then be computed using the orthogonality of Weiner Chaos expansion as
\begin{equation}\label{L2norm}
	\|f\|_{L^2(\Omega)}=\left(\sum_{m\geq 0}m!\int\left(\sum_{k\geq 0}\frac{(m+2k)!}{m!k!2^k}\int_{\mathbb{R}^{kd}}|f_{m+2k}|(t,\xi,{\bf x}_m,{\bf y}_k^{\otimes 2})d{\bf y}\right)^2 d{\bf x}\right)^{\frac{1}{2}} <\infty,
\end{equation}
The readers are referred to \cite{B-CMP-09} for more details.\\
It is easy to verify that the space defined above is dense in $L^2(\Omega\times\mathbb{R}^d)$. Denote the space consisting of all functions that satisfy condition (1) and (2) by $\tilde{M}$. In fact, any function $f(\xi)\in L^2(\Omega\times\mathbb{R}^d)$ can be written as their Wiener Chaos expansion
\begin{equation}
	f(\xi) = \sum_{m\geq 0}I_m(g_m(\xi,{\bf x}_m)).
\end{equation}
Each $g_m$ can be approximated by a function $f_m^k$, which vanishes in a set of measure at most $k^{-1}$ in the vicinity of measure $0$ set of diagonals given by the support of the distributions $\delta(x_e-x_f)$. By the change of change of coordinates in \eqref{change}, we have $g_m^{(k)}=f_m^{(k)}$ so that the It{\^o} and Stratonovich integrals agree. Define
\begin{equation}
	f^{(k)}(\xi) = \sum_{m\geq 0}\mathcal{I}_m(f_m^{(k)}).
\end{equation}
Using formula \eqref{L2norm}, we may verify that
\begin{equation}
	\|\mathcal{J}f^{(k)}\|_{L^2(\Omega)} < \infty,
\end{equation}
and
\begin{equation}
	\lim_{k\rightarrow \infty}\|f^{(k)}(\xi)-f(\xi)\|_{L^2(\Omega\times\mathbb{R}^d)} = 0.
\end{equation}
We have shown that $\tilde{M}$ is dense in $L^2(\Omega\times\mathbb{R}^d)$. Since $M$ is dense in $\tilde{M}$, it is also dense in $L^2{(\Omega\times\mathbb{R}^d)}$.
\section{General moment convergence}\label {section:general}
We now turn to the general case of \eqref{thm2}. It suffices to consider the limit
\begin{equation}
	\lim_{\varepsilon\rightarrow 0}\mathbb{E}\{\hat{u}_{\varepsilon}(t_1,\xi^{(1)})\cdots\hat{u}_{\varepsilon}(t_r,\xi^{(r)})\}.
\end{equation}
Using \eqref{ue1}, we obtain that the above expression before passing to the limit equals
\begin{equation}
\label{general}
\sum\mathcal{I}_{\varepsilon}(\boldsymbol{n}),
\end{equation}
where
\begin{equation}
	\mathcal{I}_{\varepsilon}(\boldsymbol{n}):=\mathbb{E}\{\hat{u}^{(n_1)}_{\varepsilon}(t_1,\xi^{(1)})\cdots\hat{u}^{(n_r)}_{\varepsilon}
	(t_r,\xi^{(r)})\}.
\end{equation}
The summation in \eqref{general} extends over all non-negative integer multi-indices $\boldsymbol{n}=(n_1,\cdots,n_r)$. We adopt the customary notation $|\boldsymbol{n}|=\sum_{l=1}^r n_l$. In fact, only the terms for which $|\boldsymbol{n}|$ are even do not vanish.\\
Each term appearing on the right-hand side can be written as
\begin{equation}
\begin{aligned}
	\mathcal{I}_{\varepsilon}(\boldsymbol{n}):=&(-i)^{|\boldsymbol{n}|}\exp(\sum_{l=1}^r t_{l})\sum_{\pi}\int_{\mathbb{R}^{|\boldsymbol{n}|d}} \exp\{i\sum_{l=1}^r\beta_l t_l\} \prod_{l=1}^{r} d\beta_l \int\cdots\int \prod_{l=1}^r d\bxi^{(l)}\\
	&\times\prod_{l=1}^{r}\hat{u}_0(\xi^{(l)}-\sum_{k=1}^{n_l}\xi_k^{(l)})\prod_{(ef)\in\pi}\hat{R}(\varepsilon\xi_{e})\delta(\xi_e+\xi_f)\\
	&\times\left\{\prod_{l=1}^r\prod_{k=0}^{n_l}\left[1-i\left(\left|\xi^{(l)}-\sum_{j=1}^{k} \xi_j^{(l)}\right|^{\mathfrak{m}}-\beta\right)\right]\right\}^{-1}.
\end{aligned}
\end{equation}
We denote $d\bxi^{(l)}:=d\xi_1^{(l)}\cdots d\xi_{n_l}^{(l)}$. The summation extends over all possible pairings. We can repeat the argument made in the previous section and obtain bounds of the form
\begin{equation}
	|\mathcal{I}_{\varepsilon}(\boldsymbol{n})|\leq (|\boldsymbol{n}|-1)!! \frac{(Ct^{2-\rho})^{\frac{|\boldsymbol{n}|}{2}}e^{rT})}{(\frac{|\boldsymbol{n}|}{2}!)^{2-\rho}},
\end{equation}
assuming $0\leq t_1,\dots,t_r\leq T$, where the constant $C$ does not depend on either $\varepsilon\in (0,1]$, $|\boldsymbol{n}|$ or $r$. There are ${N+r-1\choose r-1}$ non-negative integer valued multi-indices satisfying equation $|\boldsymbol{n}|=N$. We can therefore estimate
\begin{equation}
\label{generalest}
\left|\sum_{|\boldsymbol{n}|=N}\mathcal{I}_{\varepsilon}(\boldsymbol{n})\right|\leq c_N,
\end{equation}
where
\begin{equation}
c_N:={N+r-1\choose r-1}(N-1)!! \frac{(CT^{2-\rho})^{\frac{N}{2}}e^{rT}}{(\frac{N}{2}!)^{2-\rho}},
\end{equation}
and the series $\sum_{N=1}^{+\infty} c_N$ is clearly summable. We are therefore allowed to pass \eqref{general} to the limit $\varepsilon\rightarrow 0_+$. As a result, we conclude that
\begin{equation}
	\lim_{\varepsilon\rightarrow 0_+}\mathbb{E}\{ \hat{u}_{\varepsilon}(t_1,\xi^{(1)})\cdots\hat{u}_{\varepsilon}(t_r,\xi^{(r)})\}=\sum\lim_{\varepsilon\rightarrow 0_+}\mathcal{I}_{\varepsilon}(\boldsymbol{n})=\sum\bar{\mathcal{I}}(\boldsymbol{n}),
\end{equation}
where
\begin{equation}
\begin{aligned}
	\bar{\mathcal{I}}(\boldsymbol{n}):=&(-i\sigma(2\pi)^{-d})^{|\boldsymbol{n}|}\sum_{\pi}\int_{\mathbb{R}^|\boldsymbol{n}d|}
\exp\{i\sum_{l=1}^{r}\beta_l t_l\} \prod_{l=1}^r d\beta_l \int\cdots\int \prod_{l=1}^{r} d\bxi^{(l)}\\
&\times\prod_{l=1}^{r} \hat{u}_0(\xi^{(l)}-\sum_{k=1}^{n_l}\xi^{(j)}_k)\prod_{(ef)\in\pi}\delta(\xi_e+\xi_f)\\
&\times\left\{\prod_{l=1}^{r}\prod_{k=0}^{n_l}\left[1-i\left(\left|\xi^{(l)}-\sum_{j=1}^k
\xi^{(l)}_j\right|^{\mathfrak{m}}-\beta_l\right)\right]\right\}^{-1}\\
=&\,\mathbb{E}\{\hat{u}^{(n_1)}(t_1,\xi^{(1)})\cdots\hat{u}^{(n_r)}(t_r,\xi^{(r)})\}.
\end{aligned}
\end{equation}

\section{Weak convergence}
To show weak convergence, we need to first demonstrate the convergence of finite dimensional distributions, i.e., the convergence of distributions of $(\hat{u}_{\varepsilon}(t_1,\xi^{(1)}),\cdots,\hat{u}_{\varepsilon}(t_r,\xi^{(r)}))$ for an arbitrary $r\geq 1, \,t_1, \cdots, t_r \geq 0, \,\xi^{(1)},\cdots,\xi^{(r)}\in\mathbb{R}^d$. Since we have already proved the moment convergence in the previous sections, it suffices to verify the determinacy of the distributions by their moments. For simplicity we consider only the case $r=1$. It is easy to generalize to the case for arbitrary $r$. Using estimates \eqref{generalest} we obtain that
\begin{equation}
\label{mun}
|\mathbb{E}[\hat{u}(t,\xi)]^n|\leq \sum_{N=0}^{+\infty} {N+r-1\choose r-1}(N-1)!! \frac{(Ct^{2-\rho})^{\frac{N}{2}}e^{rt}}{(\frac{N}{2}!)^{2-\rho}}.
\end{equation}
Using Stirling's formula we can easily obtain that
\begin{equation}
{N+r-1 \choose r-1}\leq C(1+\frac{r-1}{N})^N (1+\frac{N}{r-1})^{r-1}.
\end{equation}
Plugging this into equation \eqref{mun} therefore gives
\begin{equation}
\begin{aligned}
|\mathbb{E}[\hat{u}(t,\xi)]^n|\leq &C 2^{r-1} e^{rt} \sum_{N=0}^{r-1} \left(1+\frac{r-1}{N}\right)^N (N-1)!! \frac{(Ct^{2-\rho})^{\frac{N}{2}}e^{rt}}{(\frac{N}{2}!)^{2-\rho}}\\
&+C 2^N e^{rt}\sum_{N=r}^{+\infty} \left(1+\frac{N}{r-1}\right)^{r-1}(N-1)!! \frac{(Ct^{2-\rho})^{\frac{N}{2}}e^{rt}}{(\frac{N}{2}!)^{2-\rho}}.
\end{aligned}
\end{equation}
The first term can be easily estimated by $Cr^r$ while the second by a constant $C$ independent of $r$. Therefore we have
\begin{equation}
	\sum_{r=1}^{+\infty} \frac{1}{\mathbb{E}[\hat{u}(t,\xi)]^{1/2r}}\geq C\sum_{n=1}^{+\infty}\frac{1}{r^{1/2}}=+\infty
\end{equation}
and the uniqueness follows from Carleman's condition. 

It remains to prove the tightness of $\{\hat{u}_\varepsilon(t,\xi)\}$ over $C([0,+\infty),\mathcal{S}'(\mathbb{R}^d))$. By \cite{Mitoma}, it suffices to prove that $\{u_{\varepsilon}^{\phi}(t):=\langle\hat{u}_{\varepsilon}(t,\cdot),\phi\rangle_{L^2(\mathbb{R}^d)},t\geq 0 \}$ for an arbitrary $\phi\in\mathcal{S}(\mathbb{R}^d)$, which, by Kolmogorov's theorem, follows from:
\begin{proposition}
For any $T>0$ and $\phi\in \mathcal{S}(\mathbb{R}^d)$ there exists a constant $C>0$ such that
\begin{equation}
\label{prop5.1}
\mathbb{E}|u_{\varepsilon}^{\phi}(t)-u_{\varepsilon}^{\phi}(s)|^2\leq C(t-s)^2, \forall \varepsilon\in(0,1],s,t\in[0,T].
\end{equation}
\end{proposition}
\begin{proof}
$\mbox{}$From equation \eqref{eq1.2}, we obtain
\begin{equation}
	u_{\varepsilon}^{\phi}(t)-u_{\varepsilon}^{\phi}(s)=i\xi^{\mathfrak{m}}\int_s^t u_{\varepsilon}^{\phi}(\tau)d\tau+\int_s^t v_{\varepsilon}^{\phi}(\tau)d\tau,
\end{equation}
and
\begin{equation}
\begin{aligned}
v_{\varepsilon}^{\phi}(\tau):=&\sum_{n\geq 0} v_{\varepsilon}^{n,\phi}(\tau),\\
v_{\varepsilon}^{n,\phi}(\tau):=&(-i)^{n+1}\varepsilon^{-\frac{d(n+1)}{2}}e^{\tau}\int_{\mathbb{R}}e^{i\beta\tau}d\beta
\int\cdots\int \prod_{k=0}^{n}\hat{q}(\xi_k)d\xi_k\int_{\mathbb{R}^d} d\xi\\
&\times \left\{\prod_{k=0}^{n}\left[1-i\left(\left|\xi-\varepsilon^{-1}\sum_{j=0}^{k}\xi_j\right|^{\mathfrak{m}}-\beta\right)\right]\right\}^{-1}
\hat{u}_0(\xi-\varepsilon^{-1}\sum_{j=0}^{n}\xi_j)\phi(\xi).
\end{aligned}
\end{equation}
As in section 4, we conclude that for all $\tau\in [0,T]$, we have $\mathbb{E}|i\xi^{\mathfrak{m}}u_{\varepsilon}^{\phi}(\tau)|^2\leq C$ and $\mathbb{E} |v_{\varepsilon}^{\phi}(\tau)|^2\leq C$ for some constant independent of $\varepsilon\in (0,1]$. Estimate \eqref{prop5.1} is a consequence of the Cauchy-Schwarz inequality.
\end{proof}
Using \eqref{un2} for any $\phi\in\mathcal{S}(\mathbb{R}^d)$, we can write
\begin{equation}
\begin{aligned}
&\langle\hat{u}^{(n)}(t,\cdot),\phi\rangle-\langle\hat{u}^{(n)}(s,\cdot),\phi\rangle\\
= &\int_s^t e^{\tau}d\tau \int_{\mathbb{R}} e^{-i\beta\tau} d\beta\int\cdots\int
\prod_{k=1}^{n}d\xi_k\int d\xi \\
&\times \left\{\prod_{k=0}^{n}\left[1-i\left(\left|\xi-\sum_{j=1}^{k}\xi_j\right|^{\mathfrak{m}}+\beta\right)\right]\right\}^{-1} \hat{u}_0(\xi-\sum_{j=1}^n\xi_j)\phi(\xi)\prod_{j=1}^{n}e^{-i\xi_j x_j} \circ\prod_{j=1}^{n}dW(x_j)\\
&-i\int_s^t e^{\tau}d\tau \int_{\mathbb{R}} \beta e^{-i\beta\tau} d\beta\int\cdots\int
\prod_{k=1}^{n}d\xi_k\int d\xi\\
&\times \left\{\prod_{k=0}^{n}\left[1-i\left(\left|\xi-\sum_{j=1}^{k}\xi_j\right|^{\mathfrak{m}}+\beta\right)\right]\right\}^{-1} \hat{u}_0(\xi-\sum_{j=1}^n\xi_j)\phi(\xi)\prod_{j=1}^{n}e^{-i\xi_j x_j} \circ\prod_{j=1}^{n}dW(x_j).\\
\end{aligned}
\end{equation}
Applying the same technique as in the proof for the $L^2(\mathbb{R}^d\times\Omega)$ setting, we check that for any $T>0$, we have
\begin{equation}
	\mathbb{E}|\langle\hat{u}^{(n)}(t,\cdot),\phi\rangle-\langle\hat{u}^{(n)}(s,\cdot),\phi\rangle|^2\leq (CT)^2(t-s)^2(n!)^{(1-\rho)}, \forall n\geq 1, s,t \in [0,T]
\end{equation}
for some constant $C>0$ independent of $n$. This in turn implies that
\begin{equation}
\mathbb{E}|\langle\hat{u}(t,\cdot),\phi\rangle-\langle\hat{u}(s,\cdot),\phi\rangle|^2\leq C(t-s)^2, \forall n \geq 1
\end{equation}
on any compact set. By the Kolmogorov-Chentsov Theorem, we have that $\langle\hat{u}(t,\cdot),\phi\rangle$ is continuous almost surely.

Weak convergence of $\{\hat{u}_{\varepsilon}(t,\xi)\}$ follows from the convergence of finite dimensional distribution and tightness. Returning to the spatial space, by the Plancherel theorem, we have $\langle u(t,\cdot),\phi\rangle = \langle \hat{u}(t,\cdot),\hat{\phi}\rangle$. Hence, the process $\{u_{\varepsilon}(t,x)\}$ converges in law over $C([0,+\infty);\mathcal{S}'(\mathbb{R}^d))$ to $\{u(t,x)\}$.

\section*{Acknowledgment}
The authors would like to thank Tomasz Komorowski for useful discussions during the preparation of the manuscript. This paper was partially funded by AFOSR Grant NSSEFF- FA9550-10-1-0194 and NSF grant DMS-1108608.

\end{document}